\documentclass[12pt]{article}
\usepackage{fancyhdr, array,calc,graphicx,url,tabularx}
\usepackage{geometry}
\usepackage{latexsym}
\usepackage{amssymb}
\usepackage{amsmath}
\usepackage{makeidx}
\usepackage{enumerate}
\usepackage{verbatim}
\usepackage{tikz}
\usepackage[colorlinks=true,linkcolor=blue]{hyperref}

\usepackage{changepage}

\makeindex

{
   \end{minipage}
   \vspace*{\stretch{3}}
   \clearpage
}

\renewcommand{\gg}{\gamma}

\newcommand{\bR}{{\mathbb{R}}}

\newcommand{\rest}{\restriction}

\newcommand{\card}[1]{{\vert #1 \vert} }

\renewcommand{\models}{\vDash}
\newcommand{\powerset}{{\wp}}
\newcommand{\dom}{{\rm dom}}
\newcommand{\rge}{{\rm rge}}
\newcommand{\cp}{{\rm crit }}

\newcommand{\cf}{{\rm cf}}

\newtheorem{theorem}{Theorem}[section]
\newtheorem{proposition}[theorem]{Proposition}
\newtheorem{definition}[theorem]{Definition}

\newtheorem{lemma}[theorem]{Lemma}
\newtheorem{corollary}[theorem]{Corollary}

\newtheorem{question}[theorem]{Question}

\numberwithin{figure}{section}

\newenvironment{proof}{{\it{
Proof.}}}{\nopagebreak\mbox{}{\hfill$\square$}
\par\bigskip}
\newenvironment{Claim}{\noindent{\it
Claim} {}}{}
\newenvironment{Case}{\noindent{\it
Case}}{}

\newcommand{\rprop}[1]{Proposition~\ref{#1}}
\newcommand{\rthm}[1]{Theorem~\ref{#1}}
\newcommand{\rlem}[1]{Lemma~\ref{#1}}

\newcommand{\rcor}[1]{Corollary~\ref{#1}}
\newcommand{\rdef}[1]{Definition~\ref{#1}}

\def\inseg{\trianglelefteq}

\def\k{\kappa}
\def\a{\alpha}
\def\b{\beta}
\def\d{\delta}

\def\l{\lambda}

\def\P{{\mathcal{P} }}
\def\W{{\mathcal{W} }}
\def\Q{{\mathcal{ Q}}}
\def\mH{{\mathcal{ H}}}
\def\K{{\mathcal{ K}}}

\def\R{{\mathcal R}}
\def\X{{\mathbb X}}
\def\H{{\rm{HOD}}}
\def\M{{\mathcal{M}}}
\def\N{{\mathcal{N}}}

\def\T {{\mathcal{T}}}
\def\U{{\mathcal{U}}}
\def\S{{\mathcal{S}}}

\def\X{{\mathcal{X}}}
\def\Y{{\mathcal{Y}}}

\def\VT{{\vec{\mathcal{T}}}}

\def\card#1{\left|#1\right|}

\def\and{\mathrel{\kern1pt\&\kern1pt}}

\def\inseg{\triangleleft}
\def\insegeq{\trianglelefteq}

\def\<#1>{\langle\,#1\,\rangle}

 \input xy
 \xyoption{all}
       
\title{Negative results on precipitous ideals on $\omega_1$
\thanks{The author's research was partially supported by the NSF Career Award DMS-1352034.}}

\author{Grigor Sargsyan\\
Institute of Mathematics\\
Polish Academy of Sciences\\
IMPAN}

\date{\today}

\begin{document}

\maketitle

\begin{abstract} We show that in many extender models, e.g. the minimal one with infinitely many Woodin cardinals or the minimal with a Woodin cardinal that is a limit of Woodin cardinals, there are no generic embeddings with critical point $\omega_1$ that resemble the stationary tower at the second Woodin cardinal. The meaning of ``resemble" is made precise in the paper (see \rdef{st-like}).
\end{abstract}

Given an ideal $\mathcal{I}$ on a cardinal $\k$, let $\mathbb{P}_\mathcal{I}=\powerset(\k)/\mathcal{I}$. Forcing with $\mathbb{P}_{\mathcal{I}}$ adds a $V$-ultrafilter on $\k$. An ideal $\mathcal{I}$ on $\k$ is called \textit{precipitous} if whenever $G\subseteq \powerset(\k)$ is a $\mathbb{P}_{\mathcal{I}}$-generic ultrafilter, $Ult(V, G)$ is well-founded.  $\mathcal{I}$ is $\l$-complete if for any $\gg<\l$ and $(A_\a: \a<\gg)\subseteq \mathcal{I}$, $\cup_{\a<\gg} A_\a\in \mathcal{I}$. If $\mathcal{I}$ is a $\l$-complete precipitous ideal on $\k$  then the generic embedding, $j: V\rightarrow Ult(V, G)$, produced by $\mathcal{I}$ has a critical point $\geq \l$.

It is mentioned in \cite{MartinMax} that Jech asked whether supercompact cardinals imply that the non-stationary ideal on some cardinal $\k$ is precipitous. Theorem 33 of \cite{MartinMax} shows that this is not the case, as any normal precipitous ideal can be destroyed in a forcing extension. However, the following question remained open. 

\begin{question} Do large cardinals imply that there exists a precipitous ideal on $\omega_1$ or on other regular cardinals?
\end{question}

It is in fact not hard to show that  sufficiently nice extender models do not carry precipitous ideals on $\omega_1$.  \rthm{woodin's theorem} was independently discovered by many inner model theorists. The proof generalizes to obtain stronger results on non-existence of precipitous ideals. To the author's best knowledge, these results are unpublished and not due to the author. Because of this we will not dwell on them and will just give the prototypical argument. 

\begin{theorem}\label{woodin's theorem} Suppose $\W$ is a countable $\omega_1+1$-iterable mouse and $\kappa$ is a successor cardinal of $\W$ such that  $\W\models ``\k^{+2}$ exists" and there is $E\in \vec{E}^\W$ such that $\cp(E)>(\k^{++})^\W$. Then $\W\models ``\k$ doesn't carry a $\k$-complete precipitous ideal".
\end{theorem}
\begin{proof} Let $\phi$ be the following sentence (in the language of premice, $\vec{E}$ is used for the extender sequence): there exists $\nu$ such that
\begin{enumerate}
\item $\nu$ is a successor cardinal,
\item $\nu^{++}$ exists,
\item there is an extender $E\in \vec{E}$ with $\cp(E)>\nu^{++}$,
\item there is a $\nu$-complete precipitious ideal on $\nu$.
\end{enumerate}

Towards a contradiction assume that in $\W$, $\mathcal{I}$ is a $\k$-complete precipitous ideal on $\k$. Let $E\in \vec{E}^\W$ be the least such that $\cp(E)>(\k^{++})^\W$ and let $\b=index^\W(E)$. Then $\mathcal{I}\in \W||\b$. It follows that without loss of generality we can assume that $\W=\W||\b$, $\W$ is sound and $\rho_\omega(\W)=\omega$ and if $\W'\inseg \W$ then $\W'\models \neg \phi$. 

 Let $j: \W\rightarrow \N\subseteq \W[g]$ be a generic embedding given by $\mathcal{I}$. Set $\l=(\k^{++})^\W$ and let $\eta$ be the predecessor of $\k$. We have that \\\\
 (1) $\l$ is a cardinal of $\W[g]$ as $\card{\powerset(\k)}^\W=(\k^{+})^\W$.\\
 (2) $j(\l)=\l$ and $j(\k)=(\eta^+)^\N$.\\
 (3) $\W|\l\not \insegeq j(\W|\l)$.\\\\
Notice that in $V$ there is a real $x$ that codes a premouse $\Q$ and an elementary embedding $\pi: \Q\rightarrow j(\W|\l)$ such that 
 \begin{enumerate}
 \item $\W|\eta\insegeq \Q$, 
 \item $\pi\rest (\eta+1)=id$,
 \item $\Q\not \insegeq j(\W|\l)$,
 \item $\pi(\l)=\l$.
 \end{enumerate}
  Take for instance a real that codes $(\W|\l, j\rest \W|\l)$. It follows by absoluteness that there is such a real in $\N^{Coll(\omega, j(\l))}$.  
  
  Let $h\subseteq Coll(\omega, \l)$ be $\W$-generic. It follows from elementarity  that there is a pair $(\R, \sigma)\in \W[h]$ such that
  \begin{enumerate}
  \item $\W|\eta\insegeq \R$,
  \item $\sigma:\R\rightarrow \W|\l$ is elementary, 
  \item $\sigma\rest (\eta+1)=id$, 
  \item $\R\not =\W|\l$, and
  \item $\sigma(\l)=\l$.
  \end{enumerate} 
Because $\W$ is $\omega_1+1$-iterable, it follows that the phalanx $(\W, \R, \eta+1)$ is $\omega_1+1$-iterable. We then compare $\W$ with $(\W, \R, \eta+1)$. Let $\T$ and $\U$ be the iteration trees on $\W$ and $(\W, \R, \eta+1)$ that this comparison process produces. Because $\rho_\omega(\W)=\omega$ and $\W$ is sound, the last model of $\U$ is on the top of $\R$. Let $\R'$ be the last model of $\U$ and $\W'$ be the last model of $\T$.

Because $\sigma: \R\rightarrow \W|\l$ is elementary, it follows that $(\eta^{++})^\R$ is the largest cardinal of $\R$. We now want to argue that $\R'=\R$. Assume not. Let $E$ be the first extender of $\U$ that is used on the $\R$-to-$\R'$ branch. Then because $\cp(E)>\eta$, we must have that there is a drop on the $\R$-to-$\R'$ branch of $\U$, and hence $\R'$ is not sound. It follows that we must have that $\W'=\W$ and $\W\insegeq \R'$ implying that $\W\insegeq \R$. But then $\sigma(\W)\inseg \W$ contradicting the minimality of $\W$.

We thus have that $\R'=\R$  and $\R\insegeq \W'$. We now want to argue that $\R\insegeq \W$. Assume that this is false. Let $F=E_0^\T$ and let $\b=index^{\W}(F)$. It follows that\\\\
(4) $\b\geq (\eta^+)^\R$\footnote{Because $\cp(\sigma)$ is a cardinal of $\R$, $\cp(\sigma)\geq (\eta^+)^\R$.} and $\b<Ord\cap \R=\l$, \\
(5) $\b$ is a cardinal of $\W'$ and \\
(6) $\W'|\b=\W|\b$. \\
(7) $\b\leq (\eta^{++})^\R$.\\\\
(7) is a consequence of (4), (5) and (6). Thus, it follows from (6) that we have two possibilities: either $\b=(\eta^{+})^\R$ or $\b=(\eta^{++})^\R$ .

Suppose now that $\b=(\eta^{++})^\R$ and set $\nu=(\eta^{+++})^{Ult(\W, F)}$. Since $\b<\l$, we have that $\nu<\l$. Since $\nu'=_{def}(\eta^{+++})^{\W'}\leq \nu$\footnote{Set $(\eta^{+++})^{\W'}=Ord\cap \W'$  in case $(\eta^{+++})^{\W'}$ is undefined.} and $(\eta^{++})^\R$ is the largest cardinal of $\R$, we have that $\R\insegeq \W'|\nu'$. Thus, $\l\leq \nu'\leq \nu<\l$, contradiction.

Assume then that $\b=(\eta^+)^\R$. Notice that $\R$ cannot be an initial segment of $Ult(\W, F)|\l$ because $\pi_F(\cp(F))\in (\eta, \l)$ is an inaccessible cardinal of $Ult(\W, F)|\l$ whereas $(\eta^{++})^\R$ is the largest cardinal of $\R$ and $\l\subseteq \R$. Thus, $G=_{def}E_1^\T$ is defined and if $\gg$ is the index of $G$ in $\M_1^\T$ then\\\\
(8) $\R|(\eta^{++})^\R\insegeq \M_1^\T|\gg$ (notice that $\M_1^\T|\gg=\W'|\gg$ and $\gg$ is a cardinal of $\W'$).\\\\
(8) now implies that $\R\insegeq \M_2^\T|(\eta^{+++})^{\M_2^\T}$, and therefore, $\l\leq (\eta^{+++})^{\M_2^\T}$. However, as it was the case with $\b$, $\gg<\l$ implying that $(\eta^{+++})^{\M_2^\T}<\l$.

We thus have shown that $\W'=\W$ and that consequently, $\R\insegeq \W$. Therefore, $\R=\W|\l$ contradicting the fact that $\R\not =\W|\l$.
\end{proof}

Woodin showed that \textit{strong condensation}, an axiom that he formulated, implies the non-existence of precipitous ideals on $\omega_1$ and cardinals below the least inaccessible cardinal. The proof is very similar to the one we gave above (see \cite[Definition 8.5]{Woodin} and \cite[Corollary 8.9]{Woodin}). The authors of \cite{SchVel} say that Steel showed that in some extender models, $\k$ carries a precipitous ideal if and only if it is measurable. The authors of \cite{SchVel} showed that in the minimal extender model with Woodin cardinal that is itself a limit of Woodin cardinals $\omega_1$ does not carry precipitous ideal (see \cite[Corollary 4]{SchVel}). The authors of \cite{ClavSch} showed that if the extender model is a model of $V=K$ then $\k$ carries a precipitous ideal if and only if it is a measurable cardinal (see \cite[Theorem 0.3]{ClavSch}).

The proof of our main theorem, \rthm{wlw neg}, uses a different type of argument that is not based on condensation. It is not clear to us how to prove \rthm{wlw neg} via condensation-like arguments or arguments based on the core model.

It is a well-known result of Woodin that if there is a Woodin cardinal $\d$ then letting $\mathbb{Q}_\d$ be the countable stationary tower forcing associated to $\d$ (see \cite{Stationarytower}), there is $G\subseteq \mathbb{Q}_\d$ and an embedding $j: V\rightarrow M \subseteq V[G]$ definable in $V[G]$ such that
\begin{enumerate}
\item $\cp(j)=\omega_1$ and
\item $V[G]\models M^\omega\subseteq M$.
\end{enumerate} 

The question on the existence of precipitous ideals on $\omega_1$ can be interpreted in at least two ways. One, of course, is the most direct interpretation. However, it can also be perceived as a question on the existence of generic embeddings that resemble the stationary tower embedding but are produced via small forcing, smaller than the size of the least Woodin cardinal. 

In this paper, we investigate this interpretation of the question.

\begin{definition}\label{st-like} Suppose $\d$ is a Woodin cardinal which is not a limit of Woodin cardinals. Let $\mu$ be the supremum of Woodin cardinals $<\d$. We say there is a stationary-tower like embedding (st-like-embedding) below $\d$ if there is a partial ordering $\mathbb{P}$ such that  whenever $g\subseteq \mathbb{P}$ is generic, 
\begin{enumerate}
\item $\mu< \card{\mathbb{P}}<\d$,
\item $(\mu^+)^V<\omega_1^{V[g]}$, 
\item in $V[g]$, there is an elementary embedding $j: V\rightarrow M\subseteq V[g]$ with the property that $\cp(j)=\omega_1$, $\mathbb{R}^{V[g]}\subseteq M$ and for some regular cardinal $\nu<\d$, 
\begin{center}
$M=\{ j(f)(s): s\in [\nu]^{<\omega}, f: [\nu]^{\card{s}}\rightarrow V$ and  $f\in V\}$\footnote{This condition says that $M$ is an ultrapower.}.
\end{center}
\end{enumerate}
\end{definition}
The last portion of clause 3 above implies that $M=Ult(V, E)$ where
\begin{center}
$E=\{(s, A): s\in j(A)\cap [\nu]^{<\omega} \wedge A\subseteq [\nu]^{\card{s}}\}$.
\end{center}
It is worth noting that $E$ may not be a short extender, and in the case of the countable stationary tower, it is not a short extender. If $\mathbb{Q}_{<\d}$ is the countable stationary tower at $\d$ and $j: V\rightarrow M$ is a generic ultrapower by some generic $G\subseteq \mathbb{Q}_{<\d}$ then setting $\nu=\d^{+}$, if $E$ is the $(\nu, \nu)$ extender derived from $j$, we have that $Ult(V, E)$ and $M$ agree on subsets of $j(\d)$\footnote{We have that for each $\a<\nu$, $\card{j(\a)}^V\leq \d$ implying that $j(\nu)=\nu$.}.

The main question we deal with in this paper is the following.

\begin{question}\label{main ques} Assume there is a Woodin cardinal $\d$. Is there an st-like-embedding below $\d$?
\end{question}

The following is the main theorem of this paper.

\begin{theorem}\label{wlw neg} Let $\M$ be the minimal mouse with a Woodin cardinal that is a limit of Woodin cardinals. Let $\d$ be the second Woodin cardinal of $\M$. Then there is no st-like-embedding below $\d$.
\end{theorem}

We will need the following proposition.

\begin{proposition}\label{reflecting strong} Suppose $\d$ is a Woodin cardinal which is not a limit of Woodin cardinals, and let  $\mu$ be the supremum of the Woodin cardinals $<\d$. Let $\k<\d$ be the least $<\d$-strong cardinal and let $\xi$ be the least such there is a poset $\mathbb{P}\in V_\xi$ witnessing that there is an st-like-embedding below $\d$. Then $\xi<\k$.
\end{proposition}
\begin{proof} Let $g\subseteq \mathbb{P}$ be generic and let $j: V\rightarrow M$ be the st-like embedding  in $V[g]$. Let $\nu_0<\d$ be such that 
\begin{center}
$M=\{ j(f)(s): s\in [\nu_0]^{<\omega}, f: [\nu_0]^{\card{s}} \rightarrow V$ and $ f\in V\}$.
\end{center}
Let $\nu\in (\max(\nu_0, \xi), \d)$ be an inaccessible cardinal and let $E$ be the $(\nu_0, \nu_0)$-extender derived from $j$. More precisely, 
\begin{center}
$E=\{(s, A): s\in[\nu_0]^{<\omega}, A\subseteq [\nu_0]^{\card{s}}$ and $s\in j(A)\}$.
\end{center}
Let $F$ be an extender with critical point $\k$ witnessing that $\k$ is $\nu$ strong. Set $W=Ult(V, F)$. We write $W_\a$ for $V_\a^W$. It follows that $E\in W[g]$ and because $Ult(V, E)$ is well-founded, $Ult(W, E)$ is also well-founded.

 It is now not hard to verify that $\pi^W_E: W\rightarrow Ult(W, E)$ is a st-like-embedding below $\d$ (in $W[g]$). Because $\mathbb{P}\in W_{\pi_F(\k)}$, we have that $\pi_F(\xi)<\pi_F(\k)$. Hence, $\xi<\k$.
\end{proof}

Upon seeing the results of this paper, Woodin informed us that he already knew that in extender models there is no st-like-embedding below the first Woodin cardinal (in fact condensation style arguments give this). He also informed us that the answer was not known for the second Woodin cardinal and beyond. We could have chosen any Woodin cardinal $\d$ such that the least cardinal $\k$ that is $<\d$-strong is not a limit of Woodin cardinals. Our proof has all the main ideas, and this is not a vanity contest. Thus, we chose to work with the second Woodin cardinal. 

We have not tried to prove results for overlapped Woodins, and believe that this is an interesting project. The methods of \cite{TPIDIMT} are probably relevant to this project.

Our methods are methods developed by inner model theorists for the last 60 years or so. We rely heavily on the writings of Mitchell and Steel. Readers familiar with the papers \cite{FSIT} and \cite{DMATM} can see their influence on the current paper. 

We started thinking about generic embeddings in extender models because of Mathew Foreman. He informed us that it is not known if large cardinals imply the existence of precipitous ideals on $\omega_1$. We thank him for asking us this question. 

Our motivation was just to show that inner model theory is a subject relevant to combinatorial set theory in a sense that a great deal of combinatorics beyond principles such as $\Diamond$ and $\square$ can be investigated and understood inside inner models. One only needs to try.

Nevertheless, we do agree with the view that the internal combinatorial structure of extender models have not been very extensively studied beyond \cite{SchZem}. However, there are several papers in print that do investigate the internal structure of mice in different ways than \cite{SchZem} does. For instance, \cite{FarmSuslin} characterizes homogeneously Suslin sets in extender models, and \cite{VarsovianI} investigates grounds of certain types of extender models.

\textbf{Acknowledgments.} The author is gratful to the referee for many useful comments. The author's work was partially supported by the NSF Career Award DMS-1352034 and by the NSF Award DMS-1954149. The author is grateful to the Institute of Mathematics of the Polish Academy of Sciences for hosting him during the academic year 2020-2021, and for providing an excellent environment for conducting research. The author finished the paper while being at the institute.  

\section{On $S$-reconstructible operators}\label{s-reconstructible operators}

Here we discuss some facts that describe the internal structure of a large class of mice.  Suppose that
\begin{enumerate}
\item $\M$ is a class size mouse over some set $x$ satisfying a sentence $\phi$,
\item there is no active level $\R\insegeq \M$ such that if $E$ is the last extender of $\R$ then $\R|\cp(E)\models \phi$ and
\item there is an active mouse $\R$ such that if $E$ is the last extender of $\R$ then $\R|\cp(E)\models \phi$.
\end{enumerate}
Clause 3 above implies that $\M$ has a club of indiscernibles. We then say $\M$ is the \textit{minimal class size} $x$-mouse satisfying $\phi$ if $\M$ is the hull of a club of indiscernibles. It is one of the most celebrated theorems of inner model theory that if there a minimal class size $x$-mouse satisfying $\phi$ then it is unique. This can be shown via a standard comparison argument (see for instance \cite[Theorem 3.11]{OIMT}). Just notice that if $\M$ and $\N$ are both minimal class size $x$-mice satisfying $\phi$ then their comparison has a club of fixed points all of which are indiscernibles.  

We say $\mathbb{M}: V\rightarrow V$ is a mouse operator if for some formula $\phi$, 
\begin{enumerate}
\item $\dom(\mathbb{M})= \{x: L_{\omega}[x]\models ``x$ is wellordered"$\}\cap \{x: $there is a minimal class size $x$-mouse satisfying $\phi\}$,
\item for each $x\in dom(\mathbb{M})$, $\mathbb{M}(x)$ is the minimal class size mouse satisfying $\phi$.
\end{enumerate}
 We also say that $\mathbb{M}$ is determined by $\phi$ and denote it by $\mathbb{M}_\phi$. When $\phi$ is clear from context we drop it from our notation, and for $x\in dom(\mathbb{M})$, we let $\M(x)=\mathbb{M}(x)$. We say $\mathbb{M}$ is total on a set $X$ if $\M(x)$ is defined for every $x\in X\cap \{x: L_{\omega}[x]\models ``x$ is wellordered"$\}$. 
 
We assume familiarity with \cite{FSIT} or with \cite{OIMT}. In particular, familiarity with \cite[Chapter 11]{FSIT} will be very helpful. Recall from \cite[Chapter 11]{FSIT} that the extenders used in the fully backgrounded construction are total and hence, have measurable critical points in the sense that if $\N$ is a model appearing in the fully backgrounded construction and $E\in \vec{E}^\N$ is a total extender then $\cp(E)$ is a measurable cardinal of $\N$. Also, Lemma 11.1, Lemma 11.2 and Theorem 11.3 of \cite{FSIT} are very important for us. When we talk about fully backgrounded construction done inside a structure with a distinguished extender sequence, we tacitly assume that all extenders come from this extender sequence.


\begin{definition}\label{s-constructible} We say $\mathbb{M}_\phi$ is an $\S$-reconstructible mouse operator  if 
\begin{enumerate}
\item $dom(\mathbb{M}_\phi)=\{a\in HC: L_\omega[a]\models ``a$ is well-ordered"$\}$,
\item for each $a\in dom(\mathbb{M}_\phi)$, $\M_\phi(a)$ has infinitely many Woodin cardinals the first $\omega$ of which are $(\delta_{a, i}:i\in \omega)$, 
\item for each $a\in dom(\mathbb{M}_\phi)$, for each $i\in \omega$, for each $\M_\phi(a)$-generic $g$ for a poset of $\M_\phi(a)$-size $<\d_{a, i}$,  for each $x\in (\M_\phi(a)|\d_{a, i})[g]$,  and for each  $\eta<\delta_{a, i}$ such that $x\in (\M_\phi(a)|\eta)[g]$, letting 
\begin{enumerate}
\item $\P$ be the output of the fully backgrounded construction of $(\M_\phi(a)|\d_{a, i})[g]$ done over $x$ using extenders with critical points greater than $\eta$ and 
\item $\N$ be the result of an $S$-construction that translates $\M_\phi(a)$ into an $x$-mouse over $\P$, 
\end{enumerate}
$\N\models \phi$. 
\end{enumerate}
\end{definition}

$S$ constructions are standard constructions in inner model theory. They were first considered by John Steel (hence the $``S"$). The first known full treatment of $S$ constructions was presented in \cite{Selfiter}, where, for some truly unfortunate though fully understandable reasons\footnote{Notice the $``S"$ in the last names of all the people involved in this bustiness.}, they were called $P$ constructions where $P$ stands for nothing in particular. The reader can also consult \cite[Chapter 3.8]{ATHM}.

Our goal is to consider two particular kinds of mice, $\M_\omega$ and $\M_{wlw}$. The first is the minimal class size mouse with $\omega$ Woodin cardinals, and the second is the minimal class size mouse with a Woodin cardinal $\d$ that is a limit of  Woodin cardinals. We will prove our theorems for $S$-reconstructible mice that have the \textit{internal covering property} (see \rdef{internal covering}). It is straightforward to check that both $\M_\omega$ and $\M_{wlw}$ satisfy our definition of $S$-reconstructible. Later we will show that they also satisfy the internal covering property (see \rthm{m has internal covering}).

Suppose $\mathbb{M}_\phi$ is an $S$-reconstructible mouse operator. Given $a\in dom(\mathbb{M}_\phi)$, we let 
\begin{center}
$\W(a)=\M_\phi(a)|(\card{a}^+)^{\M_\phi(a)}$.
\end{center}
 We think of $\W$ as a function whose domain is $dom(\mathbb{M}_\phi)$. Given a transitive set $N$, let $\W^N=\W\rest N$. The following is a corollary to our definition. In general, the results of this section are not new and reformulations of similar results that appeared in \cite{ATHM} and in \cite{DMATM} (for example see \cite[Chapter 3.1]{ATHM} and \cite[Theorem 5.1]{DMATM}).

\begin{corollary}\label{capturing w} Suppose $\mathbb{M}_\phi$ is an $S$-reconstructible mouse operator. Fix $a\in dom(\mathbb{M}_\phi)$ and $i\in \omega$, and set $\M=\M_\phi(a)$ and $\d=\d_{a, i}$. Then $\W^{\M|\d}$ is uniformly definable over $\M|\d$. Moreover, there is a formula $\psi$ with the property that for any poset $\mathbb{P}\in \M|\d$, for any $\M$-generic $g\subseteq \mathbb{P}$, for any $x\in HC^{\M[g]}$ and for any $\R$,
\begin{center}
$\R\insegeq \W(x)$ if and only if $\R\in \M|\d[g]$ and $\M|\d[g]\models \psi[x, \R]$. 
\end{center}
\end{corollary}

It is clear what $\psi$ must be, it is just the formula defining the fully backgrounded constructions. Note that the language of $\M$ has a symbol for the extender sequence of $\M$, and so $\psi$ may mention the extender sequence of $\M|\d$. Results of Schlutzenberg suggest that $\W$ maybe even be definable over the universe of $\M|\d$ (see \cite{FarmThesis}). However, we do not need such fine calculations. Below we give an example of such a $\psi$. Set \\\\
$\psi[x, \R]:$  $\R$ is an $x$-premouse and there is $\l<\d$ such that for every $\eta\in (\l, \d)$, if $\P$ is the output of the fully backgrounded construction of $\M|\d[g]$ done over $x$ using extenders with critical points $>\eta$, $\R\insegeq \P$. \\\\
It can be shown that $\psi$ witnesses \rcor{capturing w}.  
 
The next results show that $\M(a)$ in fact knows some fragments of its own strategy. The first lemma is a useful and easy lemma. We let $\d_{a, -1}=0$. We do not know the origin of this lemma but it has probably been discovered by many authors independently. 

\begin{lemma}\label{no Woodins} Suppose $\mathbb{M}_\phi$ is an $S$-reconstructible mouse operator. Fix $a\in dom(\mathbb{M}_\phi)$ and $i\in \omega$, and set $\M=\M_\phi(a)$ and $\d=\d_{a, i}$. Let $\mathbb{P}\in \M|\d$ and suppose $g\subseteq \mathbb{P}$ is $\M$-generic. Let $x\in \M|\d[g]$ and $\l\in (\max(\d_{i-1}, \card{\mathbb{P}}^\M), \d)$ be such that $x\in \M|\l[g]$. Let $\P$ be the output of the fully backgrounded construction of $\M|\d[g]$ done over $x$ using extenders with critical points $>\l$. Then $\P\models ``$there are no Woodin cardinals".
\end{lemma}
\begin{proof} Towards a contradiction, assume not. Let $\eta$ be the least Woodin cardinal of $\P$. Then $\M|\eta$ is generic over $\P$ for the extender algebra at $\eta$ that uses $\eta$-generators. We claim that\\

\textit{Claim.} $\P[\M|\eta]\models ``\eta$ is a Woodin cardinal". \\\\
\begin{proof}
To see that $\P[\M|\eta]\models ``\eta$ is a Woodin cardinal", let $f: \eta\rightarrow \eta$ be a function in $\P[\M|\eta]$. Because the forcing that adds $\M|\eta$ has the $\eta$-cc, there is $h:\eta\rightarrow \eta$ in $\P$ such that for every $\a<\eta$, $f(\a)<h(\a)$. Let $E\in \vec{E}^\P$ be an extender such that $\nu=_{def}\nu_E$ is a $\P$-cardinal such that  letting $\k=\cp(E)$,
\begin{center}
$\pi^\P_E(h)(\k)<\nu$.
\end{center}
Let $F\in \M$ be the resurrection of $E$. Let $\S$ be the model appearing in the construction producing $\P$ such that $F$ is added to $\S$. We have that no further model appearing in the construction projects below $\nu$ (as $\P|\nu$ is an initial segment of the final model of the construction). It follows that the canonical factor map $k: Ult(\P, E)\rightarrow \pi_F^\M(\P)$ has a critical point $\geq \nu$. Hence, $k(\pi_E^\P(h)(\k))=\pi^\P_E(h)(\k)$. It follows that
\begin{center}
$\pi^\M_F(f)(\k)<\nu\leq \nu_F$.
\end{center}
If $F\not \in \vec{E}^{\M|\eta}$ then for some $\M$-inaccessible $\xi \in (\nu, \eta)$, $F\rest \xi\in \vec{E}^{\M|\eta}$. It follows that $F\rest \xi$ witnesses Woodiness for $f$ in $\P[\M|\eta]$.  
\end{proof}
Notice now that $\W(\M|\eta)\in \P[\M|\eta]$ (this follows from $S$-reconstructibility). Hence, $\W(\M|\eta)\models ``\eta$ is a Woodin cardinal". It follows from \cite[Remark 12.7]{DMATM} that $\W(\M|\eta)$ is essentially the same as $\M|\a$ where $\a=(\eta^+)^\M$ in the case $\eta$ is a cutpoint and $\a=lh(F)$ where $F\in \vec{E}^\M$ is the first extender such that $\eta\in (\cp(F), lh(F))$. The point here is just that if there are partial extenders with critical point $\eta$ then they can be translated away.

Suppose now that $\eta$ is a cutpoint. It follows from the above discussion that $\W(\M|\eta)=\M|(\eta^+)^\M$ and hence, $\eta$ is a Woodin cardinal of $\M$, contradiction (there are no Woodin cardinals between $\l$ and $\d$). Next, suppose $\eta$ is not a cutpoint and let $F$ be the least extender overlapping it. Then, $Ult(\M, F)\models ``\eta$ is a Woodin cardinal". Hence, $\M|\cp(F)\models $``there are unboundedly many Woodin cardinals". Since $\eta$ is a cardinal of $\M$, $\cp(F)$ is also a cardinal implying that $\M$ has Woodin cardinals in the interval $(\l, \d)$, contradiction!
\end{proof} 

The next proposition shows that some fragments of the iteration strategy are universally Baire inside the mouse operators. The proof is very much like the proofs used in \cite[Chapter 3.1]{ATHM}. Below given an iteration tree $\T$ on a premouse $\N$, we let $C(\T)=\cup_{\a<lh(\T)}\M_\a^\T|lh(E_\a^\T)$. Usually $C(\T)$ is denoted by $\M(\T)$ which in this paper has a different meaning. 

 The propositions \rprop{everything is ub}, \rprop{capturing iterability}, \rprop{everything is ub g} and \rprop{w is ub} are all part of the standard literature. For example see \cite[Lemma 5.1]{DMATM} and \cite[Lemma 6.3]{FarmSuslin}\footnote{We thank the referee for providing these references.}.

\begin{proposition}\label{everything is ub} Suppose $\mathbb{M}_\phi$ is an $\S$-reconstructible mouse operator and $a\in dom(\mathbb{M}_\phi)$. Fix $i\in \omega$ and set $\d=\d_{a, i}$ and $\M=\M(a)$. Let $\Sigma$ be the unique iteration strategy of $\M$. 
 Suppose $\k\in (\d_{a, i-1}, \d)$ is an $\M$-cardinal and $\Lambda$ is the fragment of $\Sigma$ that acts on non-dropping trees on $\M|\k$ that are above $\d_{a, i-1}$. Then for each $j\in (i-1, \omega)$, $\Lambda\rest (\M|\d_{a, j})\in \M$ and whenever $g\subseteq Coll(\omega, \k)$ is $\M$-generic, $\Lambda\rest HC^{\M[g]}\in \M[g]$ and 
 \begin{center}
 $\M[g]\models `` \Lambda\rest HC^{\M[g]}$ is $\d_{a, j}$-uB". 
 \end{center}
 Furthermore, for every $j\in (i-1, \omega)$, whenever $h$ is $\M[g]$-generic for a poset of size $<\d_{a, j}$, $\Lambda\rest HC^{\M[g*h]}$ is the canonical extension of $\Lambda\rest HC^{\M[g]}$.
\end{proposition}
\begin{proof} The representative case is when $i=0$. When $i>0$ we need to work over $\M|\d_{a, i-1}$. Here we assume $i=0$. Also, the proof of the case when $j>0$ is very similar to the proof of the case when $j=0$. The only difference is that for $j>0$, the fully backgrounded constructions we consider must use extenders with critical points $>\d_{a, j-1}$. 

The fact that $\Lambda\rest HC^{\M[g]}\in \M[g]$ follows from \rcor{capturing w}. Indeed, notice that given a tree $\T$ on $\M|\k$ of limit length and according to $\Lambda$, $\Lambda(\T)$ is the unique branch $b$ such that $\Q(b, \T)$\footnote{For the definition of $\Q(b, \T)$ see \cite[Definition 6.11]{OIMT}.} exists and $\Q(b, \T)\insegeq \W(C(\T))$. Thus, to define $\Lambda$ in generic extensions of $\M$, it is enough to know that the function $\T\mapsto \W(C(\T))$ is definable on the domain of $\Lambda$. This follows from \rcor{capturing w}. For the rest of the argument we assume that $\k$ is a successor cardinal of $\M$. This assumption doesn't cause loss of generality, since if $\k$ is a limit cardinal then the conclusion of the proposition can be reached by using the conclusion of the proposition for $(\k^+)^\M$.

As pointed out by the referee, we could define $\Lambda\rest HC^{\M[g]}\in \M[g]$ as follows. Let $T$ be the tree of attempts to build a triple $(x, y, z, \pi)$ such that $x$ codes an iteration tree $\T$ on $\M|\k$, $y$ codes a cofinal well-founded branch $b$ of $\T$ such that $\Q(b, \T)$ exists, $z$ codes a countable $\N$ such that $\Q(b, \T)\in \N$, $\pi: \N\rightarrow \M|\d$ and inside $\N$, $\Q(b, \T)$ can be build via fully backgrounded constructions done over $C(\T)$.  
In this paper, especially in \rprop{capturing iterability}, \rprop{everything is ub g} and \rprop{w is ub},  we will need different sort of arguments, and so we present a somewhat more involved proof that exploits the idea of ``iterating to the background constructions". 

Next we show that $\Lambda$ is $\d$-uB in $\M[g]$. Our generically absolute definition of $\Lambda$ will also show the ``furthermore" clause of the proposition. Let $\l\in (\k, \d)$ be a cardinal, and let $\vec{C}=(\S_\xi, \R_\xi, F_\xi: \xi<\d)$ be the models of fully backgrounded construction of $\M|\d$ (or $\M|\d[g]$) done over $a$ in which extenders used have critical points $>\l$. Let $\R_\d$ be the output of $\vec{C}$. Thus, for $\a<\d$, $\R_\d||\a$ is defined to be $\R_{\xi_\a}||\a$ where $\xi_\a$ is the least such that for all $\zeta>\xi_\a$, $\R_\zeta||\a=\R_{\xi_\a}||\a$. We claim that\\

\textit{Claim 1.} for some $\xi$, $\R_\xi$ is an iterate of $\M|\k$ via a tree $\W$ such that $\pi^\W$ exists. \\\\
\begin{proof} Notice that as $\M|\k$ has no Woodin cardinals, if there was such a tree $\W$ then $\W\in \M$. Now towards a contradiction assume our claim is false.  We now compare $\M|\k$ with the construction  $(\S_\xi, \R_\xi, F_\xi: \xi<\d)$\footnote{Such comparison arguments were studied in \cite{FarmJohn}. But here we could also argue that the construction side doesn't move at all. The reader can consult \cite[Lemma 2.11]{ATHM} and \cite[Lemma 3.23]{trang2013}.}. We use $\Sigma_{\M|\k}$ on the $\M|\k$-side and $\Sigma$ on the $\M$-side. The comparison produces a tree $\T$ on $\M|\k$ according to $\Lambda$ with last model $\N$ and a non-dropping tree $\U$ on $\M$ according to $\Sigma$ with last model $\M_1$ such that $\pi^\U(\R_\d)\insegeq \N$. 

Indeed, if $\M|\k$-side lost then the comparison would have stopped before reaching stage $\pi^\U(\d)$, and so there would be some $\xi$ such that the second model of $\pi^\U(\vec{C})(\xi)$ was an iterate of $\M|\k$. This fact would be witnessed inside $\M_1$, and hence by elementarity our claim would be true in $\M$. 

 Because $\M|\k$ has no Woodin cardinals, we must have that $rud(\N)\models ``\pi^\U(\d)$ is not a Woodin cardinal". Because all initial segments of $\M|\k$ are $\phi$-small, we have that if $\R$ is the result of $S$-construction that translates $\M_1$ into a mouse over $\pi^\U(\R_\d)$ then $\N\insegeq \R$. However, because $\pi^\U(\d)$ is a Woodin cardinal of $\M_1$, we have that $\R\models ``\pi^\U(\d)$ is a Woodin cardinal", contradiction. 
\end{proof}

We now use branch condensation of $\Lambda$ to get a generically absolute definition of $\Lambda$\footnote{For more on branch condensation see \cite[Definition 0.19]{ATHM}.}. Let $g\subseteq Coll(\omega, \k)$ be $\M$-generic. For each $\l\in (\k, \d)$, let $\vec{C}_\l=(\S_\xi, \R_\xi, F_\xi: \xi<\d)$ be the output of the fully backgrounded construction of $\M|\d[g]$ done over $a$ in which extenders used have critical points $>\l$. Let $\xi_\l$ be such that $\R_{\xi_\l}$ is an iterate of $\M|\k$. This iteration must be according to $\Sigma_{\M|\k}$. Let $\pi_\l:\M|\k\rightarrow \R_{\xi_\l}$.

Suppose now $h$ is any $\M[g]$-generic for a poset of size $<\d$. Then given a non-dropping tree $\T\in \M|\d[g][h]$ on $\M|\k$ we say $\T$ is \textit{correct} if for all limit $\a<lh(\T)$, for some $\eta$ such that $\T\in \M|\eta[g][h]$ for all $\l\in (\eta, \d)$, there is an embedding $\sigma:\M^{\T}_\a\rightarrow \R_{\xi_\l}$ such that 
\begin{center}
$\pi_\l=\sigma\circ \pi^\T_{0, \a}$.
\end{center}
Given a correct tree $\T\in \M[g*h]$, we let $\phi[\T, b, \Q]$ be the statement that for some $\eta<\d$ for all $\l\in (\eta, \d)$
\begin{enumerate}
\item $b$ is a cofinal well-founded branch of $\T$ such that $\Q=\Q(b, \T)$ and
\item there is an embedding $\sigma: \M^\T_b\rightarrow \R_{\xi_\l}$ such that $\pi_\l=\sigma\circ \pi^\T_b$.
\end{enumerate}
Let $\psi[\T, b, \Q]$ be the statement that $\T$ is correct and $\phi[\T, b, \Q]$ holds. Notice that\\\\
(1) in $\M[g]$, whenever $\mathbb{P}$ is a poset of size $<\d$, $\mathbb{P}$ forces that for any correct tree $\T$ there is $b, \Q$ such that $\phi[\T, b, \Q]$.\footnote{Notice that if $h\subseteq \mathbb{P}$ is $\M[g]$-generic and $(\T, b, \Q(b, \T))$ is as in (1) then $(b, \Q(b, \T))\in \M[g][h]$ because a fully backgrounded constructions of $\M[g][h]$ done over $\mathcal{C}(\T)$ must reach $\Q(b, \T)$ implying that both $b, \Q(b, \T)\in \M[g]$. The proof of Claim 1 and the Dodd-Jensen property can be used to show the existence of $\sigma$.} \\\\\
The branch condensation of $\Lambda$ implies that such a pair $(b, \Q)$ must be unique. We then get that $\psi$ is a generically correct definition of $\Lambda$.\\

\textit{Claim 2.} For a club of countable $X\prec \M|(\d^+)^\M[g]$, letting $\pi_X: \N_X\rightarrow \M|(\d^+)^\M[g]$ be the transitive collapse of $X$, and letting $h\in \M[g]$ be $\N_X$-generic for a poset of size $<\pi_X^{-1}(\d)$, for any $(\T, b, \Q)\in \N_X[h]$,
\begin{center}
$\N_X[h]\models \psi[\T, b, \Q]$ if and only if $\M[g]\models \psi[\T, b, \Q]$.
\end{center}
\begin{proof} Left to right direction is easy and we leave it to the reader. For the other direction, suppose that $(\T, b, \Q)\in \N_X[h]$ and $\M[g]\models \psi[\T, b, \Q]$. First we claim that $\N_X[h]\models ``\T$ is correct". Suppose otherwise. Then there is a limit $\a<lh(\T)$ such that $\N_X[h]\models ``\T\rest \a$ is correct and $\T\rest \a+1$ is not correct". It follows from (1) that there is $c, \Q\in \N_X$ such that $c$ is not the branch of $\T\rest \a$ in $\T$ and $\N_X\models \psi[\T\rest \a, c, \Q]$. It follows that $\M[g]\models \psi[\T\rest \a, c, \Q]$ implying that $c$ is the branch of $\T\rest \a$ in $\T$. A similar argument shows that in fact $\N_X\models \psi[\T, b, \Q]$.
\end{proof}
\end{proof}

We now show that countable submodels also have universally Baire strategies. 

\begin{proposition}\label{capturing iterability} Suppose $\mathbb{M}_\phi$ is an $S$-reconstructible operator, $a\in dom(\mathbb{M}_\phi)$ and $i\in \omega$. Set $\M=_{def}\M_\phi(a)$ and $\d=\d_{a, i}$. Let $\pi: \N\rightarrow \M|(\d^+)^\M$ be a countable hull inside $\M$. Then $\M\models ``\N$ has a $\d$-uB iteration strategy that acts on trees above $\pi^{-1}(\d_{a, i-1})$". 
\end{proposition}
\begin{proof} Again, we only do the proof of the representative case $i=0$. Let $\Sigma$ be the unique iteration strategy of $\M$, and let $(\S_\xi, \R_\xi, F_\xi: \xi<\d)$ be the models of the fully backgrounded constructions of $\M|\d$ over $a$. We claim that for some $\xi<\d$ there is an embedding $\sigma: \N\rightarrow \R_\xi$. To build such an embedding, we compare $\N$ with the aforementioned construction of $\M|\d$. We use the $\pi$-pullback of $\Sigma$ to iterate $\N$. We claim that the construction side wins the comparisons.

To see this, assume not. We then get a tree $\T$ on $\N$ and a tree $\U$ on $\M|\d$ with last models $\N_1$ and $\M_1$ respectively such that $\pi^\U$ exists and $\pi^\U(\R_\d)\insegeq \N_1$. As there are no Woodin cardinals in $\R_\d$ (see \rlem{no Woodins}), $(\T\rest lh(\T)-1)\in \M_1$. It follows that there is a tree $\W\in \M$ on $\N$ such that $C(\W)=\R_\d$. It follows that $\M|(\omega_1)^\M\insegeq \N$, contradicting the fact that $\N$ is countable in $\M$. This contradiction shows that there is $\sigma:\N\rightarrow \R_\xi$ for some $\xi<\d$. The rest follows from \rprop{everything is ub}\footnote{Here we use the fact that the strategy of $\R_\xi$ is reducible to the portion of the strategy of the background universe that acts on non-dropping iteration trees. See \cite[Chapter 11]{FSIT}.}. It is not hard to show that the $\sigma$-pullback of the strategy of $\R_\xi$ induced by $\Sigma$ is $\d$-uB in $\M$. 
\end{proof}

We state, without a proof, a somewhat stronger version of \rprop{capturing iterability}. 

\begin{proposition}\label{everything is ub g} Suppose $\mathbb{M}_\phi$ is an $S$-reconstructible mouse operator, $a\in dom(\mathbb{M}_\phi)$ and $i\in \omega$. Set $\d=\d_{a, i}$ and $\M=\M_{\phi}(a)$. Let $\Sigma$ be the unique iteration strategy of $\M$. Suppose $g$ is $\M$-generic for a poset of size $<\d$, and let $\pi: \N[\bar{g}]\rightarrow \M|(\d^+)^\M[g]$ be a countable hull in $\M[g]$. Then $\M[g]\models ``\N$ has a $\d$-uB iteration strategy acting on trees that are above $\pi^{-1}(\d_{a, i-1})$".
\end{proposition}

The next lemma shows that for any $x$, proper initial segments of $\W(x)$ have universally Baire iterations strategies (in $\M_\phi(a)$). However, the function $x\rightarrow \W(x)$ cannot be universally Baire. For this we need to collapse the first strong cardinal of $\M(a)$. The reason is that $\W(x)$ is the set of all OD subsets of $x$ in the derived model of $\M(a)$ computed at $\d_{a, \omega}$, and this derived model is a model in which all sets are ordinal definable from a real. For more on this we refer the reader to \cite{DMATM}. 

\begin{proposition}\label{w is ub} Suppose $\mathbb{M}_\phi$ is an $S$-reconstructible mouse operator. Let $a\in dom(\mathbb{M}_\phi)$ and $i\in \omega$. Set $\M=\M_\phi(a)$ and $\d=\d_{a, i}$. Let $g$ be $\M$-generic for a poset of size $<\d$, $x\in \M|\d[g]\cap  dom(\W)$ and $\Q\insegeq \W(x)$ be such that $\rho_{\omega}(\Q)=\omega$. Let $\Lambda$ be the unique strategy of $\Q$. Then $\Lambda\rest HC^{\M[g]}\in \M|\d[g[$ and is $\d$-uB in $\M$ in the stronger sense that for any $\M[g]$-generic $h$, $\Lambda\rest HC^{\M[g*h]}$ is the canonical extension of $\Lambda\rest HC^{\M[g]}$. 
\end{proposition}
\begin{proof} We again do the proof in the representative case of $i=0$. To prove the claim fix $g, x, \Q$ as in the statement of the proposition. Let $\P$ be the output of the fully backgrounded construction of $\M|\d[g]$ done over $x$ using extenders with critical points $>\l$ where $\l$ is some cardinal $<\d$ bigger than the size of the poset.  We have that $\Q\insegeq \P$. Thus, again, the iterability of $\Q$ reduces to the iterability of some $\M|\k$ for non-dropping trees  where $\k>\l$ is a regular cardinal of $\M$. The rest of the claim follows from \rprop{everything is ub} and \rprop{capturing iterability}.
\end{proof}

\section{The Internal Covering Property}

We will need to deal with $S$-reconstructible operators with a stronger property. Recall $\W(x)$ function given by $\W(x)=\M(x)|(\card{x}^+)^{\M(x)}$.

\begin{definition}\label{internal covering} Suppose $\mathbb{M}_\phi$ is an $S$-reconstructible mouse operator. We say $\mathbb{M}_\phi$ has the internal covering property if for any $a\in dom(\mathbb{M}_\phi)$ and $i\in \omega$, letting $\M=\M(a)$ and $\d=\d_{a, i}$, for any $\mathbb{P}\in \M|\d$, $\M$-generic $g\subseteq \mathbb{P}$, $x\in \M|\d[g]$, and $\l\in (\d_{a, i-1}, \d)$ such that 
\begin{enumerate}
\item $a\in L_\omega[x]$,
\item $L_{\omega}[x]\models ``x$ is well-ordered",
\item $x\in \M|\l[g]$,
\item $\M|\d_{a, i-1}$ is generic over $\W(x)$,
\item $\mathbb{P}\in \M|\l$.
\end{enumerate} 
letting $\P$ be the output of the fully backgrounded construction of $\M|\d[g]$ done over $x$ using extenders with critical points greater than $\l$, for unboundedly many $\k<\d$, $(\k^+)^\P=(\k^+)^\M$.
\end{definition} 

Let $\mathbb{M}$ be either $x\rightarrow \M_\omega(x)$ or $x\rightarrow \M_{wlw}(x)$. Both of these operators are $S$-reconstructible. Here we show that they also have the internal covering property. 

Suppose $\M$ is a mouse, $\k$ is an inaccessible cardinal of $\M$ such that $\rho_\omega(\M)\geq \kappa$ and $\T$ is an iteration tree on $\M|\k$. We then let $\T^\M$ be the iteration tree on $\M$ that has the same tree structure as $\T$ and uses the same extenders as $\T$\footnote{We ignore issues involving ill-founded models as in cases where this notation is relevant the models are iterable.}. Similarly, given an iteration tree on $\M$ that is below $\k$, we let $\T\rest (\M|\k)$ be the iteration tree on $\M|\k$ that has the same tree structure and uses the same extenders as $\T$.  

\begin{theorem}\label{m has internal covering} $\mathbb{M}$ has the internal covering property. 
\end{theorem}
\begin{proof} 
We show that $\M=_{def}\M(\emptyset)$ satisfies the internal covering property. Here the representative case is $i=1$, so we assume $i=1$. Let $\d_0=\d_{\emptyset, 0}$ and $\d=\d_{\emptyset, 1}$. Let $\phi$ be the obvious defining formula of $\mathbb{M}$\footnote{For instance, in the case of $\M_\omega(x)$, $\phi$ is ``there are infinitely many Woodin cardinals".}.

Let $\xi$ be the sup of the Woodin cardinals of $\M$ and $g\subseteq Coll(\omega, <\xi)$ be generic over $\M$. Let $W$ be the derived model of $\M$ as computed in $\M[g]$. More precisely, $W=L(\Gamma, \bR^*)$ where $\bR^*=\cup_{\k<\xi}\bR^{\M[g\cap Coll(\omega, <\k)]}$ and $\Gamma$ is the collection of all those sets of reals $A$ of $\M(\bR^*)$ such that $L(A, \bR^*)\models AD^+$. Woodin's celebrated \textit{derived model theorem} says that $L(\Gamma, \bR^*)\models AD^+$ and in $\M(\bR^*)$, $\powerset(\bR^*)\cap W=\Gamma$. In the case of $\M=\M_\omega$, $W$ is just $L(\bR^*)$ (see \cite{DMT}). 

Working in $W$, let $\nu$ be the supremum of $OD^W$ prewellorderings of $\bR$. Below we collect some facts that can be proved using $\H$-analysis done inside $W$. The reader should consult \cite{HODCore}. Giving the complete proofs of these facts is beyond this paper. Let $\mH=\H^W_{\M|\d_0}$ and let $\Sigma$ be the unique iteration strategy of $\M$.
\begin{enumerate}
\item $V_\nu^\mH$ can be represented\footnote{We say ``represented" rather than ``is" because $\M$ is a structure in a different language. In particular, $\M|\d$ has the extender sequence of $\M|\d$ as a predicate.}  as a $\Sigma$-iterate of $\M|\d$ via an iteration that is above $\d_0$.
\item $\mH\models``\nu$ is a Woodin cardinal".
\item Suppose $\mathbb{P}\in \mH|\nu$ is a poset and $g\subseteq \mathbb{P}$ is $\mH$-generic. Then $\mH[g]\models ``\mH|\nu$ is $\nu+1$-iterable\footnote{here iterability refers to iterability with respect to the extender predicate of $\mH$.} for trees that are in $L[\mH|\nu][g]$".
\item Let $\S$ be the iterate of $\M$ such that $\mH|\nu\insegeq \S$ and if $i: \M\rightarrow \S$ is the iteration embedding then the generators of $i$ are contained inside $\nu$. Then the aforementioned strategy of $\mH|\nu$ is $\Sigma_\S\rest (L[\mH|\nu][g])$.
\end{enumerate}  
Let now $\S$ be as in clause 4 above. We want to prove now that $\S$ satisfies the internal covering. We have that $\S|\nu=\mH|\nu$.

Fix $\l\in (\d_0, \nu)$ and let $\mathbb{P}\in \S|\l$ be a poset. Let $g\subseteq \mathbb{P}$ be $\S$-generic and $x\in \M|\l[g]$ be such that $\M|\d_0$ is generic over $\W(x)$ and $L_\omega[x]\models ``x$ is well-ordered". Let $\P$ be the output of the fully backgrounded construction of $\M|\d[g]$ done over $x$ using extenders with critical points greater than $\l$. Let $\W$ be the output of the fully backgrounded construction of $\P[\M|\d_0]$ done over $\M|\d_0$ in which extenders used have critical points $>\l$. It is enough to show that in $\S$, $\W$ computes unboundedly many successors correctly.

We now compare $\mH|\nu$ with $\W$. On $\mH|\nu$ side we use the $\nu+1$-strategy in $\mH$ that acts on iteration trees in $L[\mH|\nu][g]$. Let $\Lambda$ be this strategy (which is a fragment of $\Sigma_{\S|\nu}$). Notice that $\Lambda$ induces a strategy for $\W$ (via the resurrection process described in \cite[Chapter 12]{FSIT}). Let then $\Psi$ be the strategy of $\W$ induced by $\Lambda$. Both $\Lambda$ and $\Psi$ act on trees of length $\leq \nu$ that are in $L[\mH|\nu][g]$.

The aforementioned comparison process lasts at most $\nu+1$ steps\footnote{This is a standard comparison argument for weasels. For example see Theorem 2.10 of \cite{BenJohn}.} . Suppose first that the comparison process stops in $<\nu$-steps. Let $\T$ and $\U$ be the trees on $\mH|\nu$ and $\W$ respectively  with last models $\mH_1$ and $\W_1$ respectively. We must have that both $\pi^\T$ and $\pi^\U$ exist. It follows that there is a club of $\xi$ such that $\pi^\T(\xi)=\xi=\pi^\U(\xi)$. For any such $\xi$ we have that
\begin{center} $(\xi^+)^{\mH}=(\xi^+)^{\mH_1}=(\xi^+)^{\W_1}=(\xi^+)^\W$,
\end{center}
which is what we wanted to show.

Suppose next that the comparison process lasts $\nu$-steps. Let $\T$ and $\U$ be the trees on $\mH|\nu$ and $\W$ respectively. In order to apply $\Lambda$ and $\Psi$ we must first show that both $\T, \U\in L[\mH|\nu][g]$. Notice that for any limit $\a<\nu$, the branch chosen by $\T$ and $\U$ at stage $\a$ is determined by the corresponding $\Q$-structures. More precisely, if $b=\Lambda(\T\rest \a)$ and $c=\Psi(\U\rest \a)$ then both $\Q(b, \T\rest \a)$ and $\Q(c, \U\rest \a)$ exist and are equal to respectively $\Q(\T\rest \a)$ and $\Q(\U\rest \a)$. However, since $C(\T)=C(\U)$, we must have that $\Q(\T)=\Q(\U)$. Thus, the comparison process is definable over $\mH|\nu$ (for instance see \rcor{everything is ub}). It follows that indeed $\T, \U\in L[\mH|\nu][g]$.

Let $\N$ be the output of the $S$-construction that translates $\S$ into a mouse over $\W$. We have that $\N\models \phi$. Set $b=\Lambda(\T)$ and $c=\Psi(\U)$. 

Notice that $\N$ has an iteration strategy induced by  $\Sigma_\S$, and if $\Psi^+$ is this strategy then $\Psi^+(\U^\N)=c$. Similarly, $\Sigma_\S(\T^\S)=b$.  Let $\T^+=\T^\S$ and $\U^+=\U^\N$. \\

\textit{Claim 1.}  $\M^\T_b=\M^\U_c$. \\\\
\begin{proof}
Because $\nu$ is inaccessible in $\mH$, we have that either $\M^\T_b\inseg \M^\U_c$ or $\M^\U_c\inseg \M^\T_b$\footnote{This follows from the usual comparison argument for weasels. For example see Theorem 2.10 of \cite{BenJohn}.}. Because both cases are symmetric let us deal with the case $\M^\U_c\inseg  \M^\T_b$ and leave the other case (which actually is easier as $\T$ is a tree on the universe itself) to the reader.

Again, the usual comparison argument for weasels (see for example Theorem 2.10 of \cite{BenJohn}) implies that either $\pi^{\T^+}_b(\nu)>\nu$ or there is a drop on $b$. To see this, assume that $\pi^{\T^+}_b$ exists and $\pi^{\T^+}_b(\nu)=\nu$. As $\M^\U_c\inseg \M^\T_b$, we have that $Ord\cap \M^\U_c< \nu$ which implies that there is a drop in $c$. But then $\W$ side cannot lose the comparison (once again see Theorem 2.10 of \cite{BenJohn}).

We thus have that either $\pi^{\T^+}_b(\nu)>\nu$ or there is a drop on $b$. In both cases, $\Q(b, \T)$ is defined. Now the usual comparison argument for weasels implies that $c$ doesn't have a drop and $\pi^\U_c(\nu)=\nu$. It follows that $\M^{\U^+}_c\models ``\nu$ is a Woodin cardinal". It then also follows that $\Q(b, \T)$ cannot have extenders overlapping $\nu$ as otherwise there will be Woodin cardinals in $\M^\U_c$. Hence, $\Q(b, \T)\insegeq \M(\M^\U_c)$. But since $\N\models \phi$, we have that $\M^{\U^+}_c\models \phi$, implying that in fact $\Q(b, \T)\insegeq \M^{\U^+}_c$. Hence, $\M^{\U^+}_c\models ``\nu$ is not a Woodin cardinal", contradiction.
\end{proof}

We thus assume that $\M^\T_b=\M^\U_c$. In fact, the proof of the claim above shows that both $\Q(b, \T)$ and $\Q(c, \U)$ do not exist. Hence, both $\pi^{\T^+}_b$ and $\pi^{\U^+}_c$ are defined and $\pi^{\T^+}_b(\nu)=\nu=\pi^{\U^+}_c(\nu)$.  We then have that there is, in $\mH[g]$, an $\omega$-club $C$ of $\xi<\nu$ such that
\begin{enumerate}
\item $\pi^\T_b(\xi)=\xi$ and
\item $\pi^{\U}_c[\xi]\subseteq \xi$.
\end{enumerate}

\textit{Claim 2.} For each $\xi\in C$, $\pi^\U_c(\xi)=\xi$.\\\\
\begin{proof} The claim follows from the following subclaim.\\

\textit{Subclaim.} For each $\xi\in (\l, \nu)$ such that $\cf^{\mH[g]}(\xi)=\omega$, $\cf^\W(\xi)$ is not a measurable cardinal of $\W$ whose measurability is witnessed by an extender on the sequence of $\W$. \\\\
\begin{proof} To see this, fix $\xi$ as above and suppose $\cf^\W(\xi)=\mu$ and there is a total extender $E\in \vec{E}^\W$ such that $\cp(E)=\mu$. Let $(\K_\a, \K'_\a, F_\a:\a<\nu)$ be the models of the fully backgrounded construction producing $\W$. Let $\a$ be such that for all $\b\in [\a, \nu)$, $\K_\a'|(\mu^+)^{\K_\a'}=\W|(\mu^+)^\W$. It then follows that $E$, as it is total over $\W$, has been added after stage $\a$ and therefore, $\mu$ is an inaccessible cardinal of $\mH$. But as $\cf^{\mH[g]}(\xi)=\omega$,  $\cf^{\mH[g]}(\mu)=\omega$. As $\mu>\l$, we get a contradiction!
\end{proof}

 Thus, since for $\xi\in C$, $\cf^\W(\xi)$ is not a measurable cardinal of $\W$, $\sup(\pi^\U_c[\xi])=\pi^\U_c(\xi)$ implying that for $\xi\in C$, $\pi^\U_c(\xi)=\xi$. 
\end{proof}

We now have that for each $\xi \in C$,
\begin{center} $(\xi^+)^{\mH}=(\xi^+)^{\M^\T_b}=(\xi^+)^{\M^\U_c}\geq (\xi^+)^\W$.
\end{center}
However, since $\pi^\U_c$ is continuous at each $(\xi^+)^\W$, we must indeed have that 
\begin{center} $(\xi^+)^{\mH}=(\xi^+)^{\M^\T_b}=(\xi^+)^{\M^\U_c}=(\xi^+)^\W$.
\end{center}
which is what we wanted to prove\footnote{We could have also argued using a pressing down argument instead of the argument given in the proof of Claim 2.}.
\end{proof}

The above proof shows that in fact covering holds on a stationary set. The referee has pointed out that the proof of \rthm{m has internal covering} uses arguments similar to those appearing in Chapter 3 of \cite{CMIP}.

\section{No towers resembling the stationary tower} 

\begin{theorem}\label{no strong stat tower} Suppose $\mathbb{M}_\phi$ is an $S$-reconstructible mouse operator with the internal covering property. Let $\a\in dom(\mathbb{M})$ and set $\M=\M(a)$ and $\d=\d_{a, 1}$. Then there is no st-like-embedding below $\d$.
\end{theorem}
\begin{proof} Towards a contradiction suppose $\mathbb{P}\in \M|\d$ is such that $\d_{a, 0}<\card{\mathbb{P}}<\d$ and whenever $g\subseteq \mathbb{P}$ is generic, there is an elementary embedding $j: \M\rightarrow \N\subseteq \M[g]$ in $\M[g]$ with the property that
\begin{enumerate}
\item $\cp(j)=\omega_1^\M$,
\item $\mathbb{R}^{\M[g]}\subseteq \N$,
\item $\card{(\d_{a, 0}^+)^\M}^{\M[g]}=\omega$,
\item  for some $\M$-regular cardinal $\nu_0<\d$, 
\begin{center}
$\N=\{ j(f)(s): s\in [\nu_0]^{<\omega}, f: [\nu_0]^{\card{s}}\rightarrow \M$ and  $f\in \M\}$.
\end{center}
\end{enumerate}

 Fix such a tuple $(g, \N, j)$. Let $\k=\omega_1^\M$ and let $\l=\omega_1^{\M[g]}$. Let $\Sigma$ be the unique iteration strategy of $\M$. Recall that we have set $\W(x)=\M(x)|(\card{x}^+)^{\M(x)}$.
We assume that $\mathbb{P}$ has the smallest possible rank.

Let $\R^*\insegeq \N$ be the least such that $\rho_\omega(\R^*)=\omega$ and $\M|\k\inseg \R^*$. Let $\Phi$ be the $\d$-strategy of $\R^*$ in $\N$. \rcor{capturing w} implies that $\W\rest \M|\d$ is definable over $\M|\d$. It then makes sense to write $\W^\N$ for the function given by the same definition over $\N|\d$.\\

\textit{Claim 1.} Suppose $x\in HC^{\M|\d[g]}$ is a transitive set such that $L_{\omega}[x]\models ``x$ is well-ordered". Then 
\begin{center}
$\W(x)\insegeq \W^\N(x)$
\end{center}
\begin{proof}
Let $\W^*\insegeq \W(x)$ be such that $\rho_\omega(\W^*)=Ord\cap x$. The proof of \rprop{w is ub} shows that for some $\tau\in (\card{\mathbb{P}}, \d)$, $\W^*$ appears as a model of the fully backgrounded construction of $\M|(\tau^+)^\M$ done using extenders with critical points $>\card{\mathbb{P}}$. Working in $\M[g]$, let $\pi: \bar{\M}[\bar{g}]\rightarrow (\M|(\d^+)^\M[g])$ be such that $x\in \bar{\M}[\bar{g}]$, $\tau\in rng(\pi)$ and $\bar{\M}$ is countable in $\M[g]$. Let $\bar{\tau}=\pi^{-1}(\tau)$. We then have that\\\\
(*)  the iterability of $\W^*$ reduces to the iterability of $\bar{\M}|(\bar{\tau}^+)^{\bar{\M}}$ for non-dropping trees that are above  $\card{\pi^{-1}(\mathbb{P})}^{\bar{\M}}$.\\\\
 By absoluteness we have $\sigma:\bar{\M}|(\bar{\tau}^+)^{\bar{\M}}\rightarrow j(\M|(\tau^+)^\M)$ in $\N$. It follows from \rprop{everything is ub} that $j(\M|(\tau^+)^\M)$ is $\d+1$-iterable in $\N$ (for non-dropping iterations that are above $\d_{a, 0}$), and hence $\W^*$ is $\d$-iterable in $\N$ (see (*)). Therefore, $\W^*\insegeq \W^\N(x)$\footnote{This follows from the universality of the fully backgrounded constructions, for example see \cite[Lemma 2.12]{ATHM}.}. 
\end{proof}

The proof of Claim 1 is a prototypical argument that we will use again below. Recall that given an iteration tree $\T$ on a premouse $\N$, we let $C(\T)=\cup_{\a<lh(\T)}\M_\a^\T|lh(E_\a^\T)$.\\

\textit{Claim 2.} There is a premouse $\X\in HC^{\M[g]}$ such that 
\begin{enumerate}
\item $\M|(\d_{a, 0}^+)^\M$ is generic over $\X$,
\item there is an iteration tree $\K$ on $\M|(\d_{a, 0}^+)^\M$ such that $\K$ is according to $\Sigma$ and either
\begin{enumerate}
\item the iteration embedding $\pi^\K$ exists and $\X$ is the last model of $\K$, or
\item  $\K$ is of limit length, $\W(C(\K))\models ``\d(\K)$ is a Woodin cardinal" and  $\X=\W(C(\K))$,
\end{enumerate}
\item there is a sound $\X$-premouse $\R\in \M[g]$ such that $Ord\cap \X$ is a cardinal of $\R$ and $\M[g]\models ``\R$ is not $\d$-iterable above $Ord\cap \X$",\\\\
letting $\nu$ be the Woodin cardinal of $\X$,
\item $\R\models ``\nu$ is a Woodin cardinal",
\item $rud(\R)\models ``Ord\cap \X$\footnote{This essentially says that $\nu^+$ of $\X$ is definably not a cardinal over $\R$.} is not a cardinal" and
\item $\X$ and $\R$ are countable in $\M[g]$. 
\end{enumerate} 
\begin{proof}
Let $\Y=\M|(\d_{a, 0}^+)^\M$. Working inside $\M[g]$, we compare $\Y$ with $\R^*$. $\Y$ is not $\d$-iterable inside $\M[g]$.\footnote{To see this, notice that if $\nu$ is the least $<\d_{a, 1}$-strong cardinal of $\M$, then because $\mathbb{P}$ was chosen to be of smallest rank, $\mathbb{P}\in \M|\nu$ (see \rprop{reflecting strong}). Then, there is an iteration tree $\U$ on $\Y$ of limit length such that $\d(\U)=(\nu^+)^\M$. $\U$ is the iteration to make $\M|\nu$-generic. $\U$ cannot have a branch in either $\M$ or $\M[g]$. For more on this kind of arguments see, for example, \cite{VarsovianI} or \cite[Chapter 5.1]{DirectedSystems}.} However, it follows from \rcor{capturing w} that the fragment of $\Sigma_\Y\rest (\M[g])$ that acts on \textit{short trees}, i.e. trees $\T$ for which $\Q(\T)$ exists, is in $\M[g]$. We then want to use the aforementioned fragment of $\Sigma_\Y$ for the comparison that we would like to perform. Finally, we would like to incorporate $\Y$-genericity iteration into the above mentioned comparison. More precisely, working inside $\M[g]$ we first iterate the least $\Y$-measurable cardinal $\d_{a, 0}+1$-times and get $\Y_1$ and then construct iteration trees $\T$ and $\U$ on $\Y_1$ and $\R^*$ respectively such that
\begin{enumerate}
\item $\T$ is according to the short fragment of $\Sigma_\Y\rest (\M[g])$,
\item $\U$ is according to $\Phi$ (recall that our hypothesis is that $\bR^{\M[g]}=\bR^\N$),
\item for $\a<\l$, given $\T\rest \a+1$ and $\U\rest \a+1$ we let 
\begin{enumerate}
\item $E^\T_{\a, 0}$ be the least extender, if it exists, on  the sequence of $\M^\T_\a$ that violates an identity in the relevant extender algebra,
\item $E^\T_{\a, 1}$ be the least extender, if it exists, that causes disagreement between $\M^\T_\a$ and $\M^\U_\a$,
\end{enumerate}
\item $E^\T_\a$ is defined if either $E^\T_{\a, 0}$ or $E^\T_{\a, 1}$ is defined and 
\begin{center}
$E^\T_\a=\begin{cases}
E^\T_{\a, 0}&: lh(E^\T_{\a, 0})\leq lh(E^\T_{\a, 1})\\
E^\T_{\a, 1} &: lh(E^\T_{\a, 1})\leq  lh(E^\T_{\a, 0})
\end{cases}$\end{center}
\item if $E^\T_\a=E^\T_{\a, 1}$ then $E^\U_\a=\vec{E}^{\M^\U_\a}(lh(E_\a^\T))$ and otherwise $E^\U_\a=E^\T_\a$. 
\end{enumerate}
Because $\R^*\not \insegeq \M$ we must have that $\R^*$-side must win any successful coiteration with $\Y$.  Notice that Claim 1 implies that the construction of $(\T, \U)$ can be carried out inside $\N$. It follows that our construction of $\T$ and $\U$ can last at most $\omega_1^\N$ many steps producing trees $\T$ and $\U$. If now $\U$ is of limit length then because $\Phi$ acts on $\U$, we can let $c=\Phi(\U)$. Let then $\R_1$ be either the last model of $\U$ (in case it has a last model) or $\M_c^\U$.

We then must have one of the following cases (in $\M[g]$): either
\begin{enumerate}
\item $\T$ has a last model $\Y_2$, $\pi^\T$ exists, $\Y$ is generic over the extender algebra of $\Y_2$ at $\pi^\T(\d_{a, 0})$ and $\Y_2\inseg \R_1$, or
\item $\T$ is of limit length, $\Q(\T)$ does not exist, and letting $\Y_2=\W(C(\T))$, $\Y_2\insegeq \R_1$.
\end{enumerate}

Set $\X=\Y_2$. Notice that because $\rho_\omega(\R^*)=\omega$, we must have that $rud(\R_1)\models``\zeta$ is not a Woodin cardinal" where $\zeta$ is the Woodin of $\X$. Let $\R_2\insegeq \R_1$ be the longest such that $\R_2\models ``\zeta$ is a Woodin cardinal". Let $\R\insegeq \R_2$ be the longest such that $\R\models ``Ord\cap \X$ is a cardinal". We claim that $(\X, \K, \R)$ is as desired where $\K$ is the iteration tree on  $\M|(\d_{a, 0}^+)^\M$ producing $\X$. It follows from our construction that it is enough to show that\\\\
(a) $\R$ is not $\d$-iterable in $\M[g]$ above $Ord\cap \X$, and\\
(b) $\X$ and $\R$ are countable in $\N$.\\\\
 Assume that (a) fails.  
Let $\eta\in (\d_{a, 0}, \d)$ be an $\M$-cardinal such that $\mathbb{P}\in \M|\eta$. Let $\P$ be the output of  the fully backgrounded construction of $\M|\d[g]$ done over $\X$ using extenders with critical point $>\eta$. Then $\R\not \insegeq \P$ which means that $\R$ must outiterate $\P$ (this can be shown by considering the comparison of $\R$ with the construction producing $\P$). We then have a tree $\U'$ on $\R$ such that $C(\U')=\P$ implying that $\P$ cannot compute unboundedly many successors correctly contradicting \rthm{m has internal covering}\footnote{For example, fix large enough inaccessible cardinal $\zeta$ such that $(\zeta^+)^\P=(\zeta^+)^\M$. Let $\a$ be the least such that $\P|\zeta\insegeq \M_\a^{\U'}$ (if $\zeta$ is chosen large enough, we must have that $\zeta=\a$ and $\zeta$ is a cardinal of $\M[g]$). Let $\b$ be the index of $E_\a^{\U'}$. We have that $\b>\zeta$ and $\b$ is a cardinal of $C(\U')$ and therefore, $(\zeta^+)^\M<\b$ (as $\b$ is not a regular cardinal of $\M_\a^{\U'}$). Thus, we must have that $\P|(\zeta^+)^\P\insegeq \M_\a^{\U'}$ and therefore, for some $\a'\in [0, \a)_{\U'}$, $(\zeta^+)^\M\in \rge(\pi^{\U'}_{\a', \a})$. As iteration embeddings are continuous at successor cardinals, we have that $(\zeta^+)^\M$ is not a regular cardinal in $\M[g]$, contradiction.}.

 Finally we need to show that $\X$ and $\R$ are countable in $\N$. Assume not. We then must have that the construction of $\T$ and $\U$ lasts $\omega_1^\N$ steps. We now claim that we must also have a branch for $\T$ in $\N$. Indeed, let $\pi: H\rightarrow \N|(\omega_2^\N)$ be countable in $\N$ such that $\T, \U, c \in rng(\pi)$. Then $\pi^{-1}(c)\in H$ is the branch of $\pi^{-1}(\U)$. Let $\xi=\omega_1^H$. Notice that $\Q(\T\rest \xi)$ exists and the branch of $\T\rest \xi$ chosen in $\T$ for $\T\rest \xi$ is the unique branch $d$ such that $\Q(d, \T\rest \xi)$ exists and is equal to $\Q(\T\rest \xi)$. But we have that $\pi^{-1}(\R_2)=\Q(\T\rest \xi)$, and because $\pi^{-1}(\R_2)\in H$, the branch of $\T\rest \xi$ is in $H$. Let $d$ be this branch. It then follows that $\pi(d)$ is a branch of $\T$. The usual comparison argument now implies that the iteration must have lasted $<\omega_1^\N$ steps.
\end{proof}

Let $\Lambda$ be the strategy of $\R$ in $\N$. We have that $\Lambda$ is a $(\d, \d)$-strategy. We would like to find a $\Lambda$-iterate $\S$ of $\R$ such that $\S$ is a \textit{minimal counterexample} to $\d$-iterability. Below we define what this notion means. 

Given a finite stack of normal trees $\VT\in \M|\d[g]$ on $\R$, we say $\VT$ is $\Lambda$-correct if in $\M[g]$, there is a club of countable  $X\prec \M|(\d^+)^\M[g]$ such that letting $\pi_X: N_X\rightarrow \M|(\d^+)^\M[g]$ be the transitive collapse, $\pi^{-1}_X(\VT)$ is according to $\Lambda$. We now look for an iterate of $\R$ that is a \textit{minimal counterexample} to $<\d$-iterability among $\Lambda$-correct iterates of $\R$. Below we make the notion more precise.

Suppose $\VT\in \M|\d[g]$ is a finite $\Lambda$-correct stack on $\R$ with last model $\K$. Let $\S\insegeq \K$. We say $(\VT, \S)$ is a minimal counterexample to $\d$-iterability if there is an $\S$-cardinal $\eta$ such that
\begin{enumerate}
\item  $\eta$ is a strong cutpoint in $\S$\footnote{This means that for all $E\in \vec{E}^\S$, if $lh(E)\geq \eta$ then $\cp(E)>\eta$.},
\item $\rho_\omega(\S)\leq \eta$ and $\S$ is $\eta$-sound,
\item in $\M|\d[g]$, $\S$ is not $\d$-iterable above $\eta$,
\item whenever $\U\in \M|\d[g]$ is a normal tree on $\S$ above $\eta$ with last model $\W^*$ such that $\VT^\frown \U$ is $\Lambda$-correct, for any $\W$ such that $\S|\eta\inseg \W\inseg \W^*$ and for any $\W$-cardinal $\nu$ such that $\nu$ is a strong cutpoint of $\W$ and $\rho_\omega(\W)\leq \nu$, $\M[g]\models ``\W$ is $\d$-iterable above $\nu$".
\end{enumerate}
It is not difficult to see that there is a minimal counterexample to $\d$-iterability. Towards a contradiction, assume there is no minimal counterexample to $\d$-iterability. We know that $\R$ is not a minimal counterexample to $\d$-iterability. We can then construct a sequence $(\R_i^+, \R_i, \T_i, \nu_i: i\in [1, \omega))$ such that\footnote{As pointed out by the referee, arguments similar to this one have appeared in the literature. For example, \cite[Lemma 1.8]{Selfiter} is similar to our argument.}.  
\begin{enumerate}
\item $\R_i$ is a $\nu_i$-sound mouse over $\R_i|\nu_i$ such that $\rho_\omega(\R_i)\leq \nu_i$,
\item $\nu_i$ is a strong cutpoint of $\R_i$,
\item $\T_i$ is a tree on $\R_i$ above $\nu_i$ such that $\oplus_{k\leq i}\T_k$ is $\Lambda$-correct,
\item $\R_{i+1}^+$ is the last model of $\T_i$,
\item $\R_{i+1}\inseg \R_{i+1}^+$ is such that for some $\nu_{i+1}$, $\R_{i+1}$ is a $\nu_{i+1}$-sound mouse over $\R_{i+1}|\nu_{i+1}$ such that $\rho_\omega(\R_{i+1})\leq \nu_{i+1}$ and $\R_{i+1}$ is not $\d$-iterable above $\nu_{i+1}$ in $\M[g]$.
\end{enumerate}
Suppose then $X\prec\M|(\d^+)^\M[g]$ is such that it witnesses that for each $i$, $\oplus_{k\leq i}\T_k$ is $\Lambda$-correct. It follows that $\pi_X^{-1}(\oplus_{i\in \omega}\T_i)$ witnesses that $\Lambda$ is not an iteration strategy for $\R$.

Let now $\S$ be a minimal counterexample to $\d$-iterability and let $\VT$ be the finite $\Lambda$-correct stack on $\R$ producing $\S$. Thus, $\S$ is an initial segment of the last model of $\VT$. Let $\eta$ be an $\S$-cardinal witnessing that $\S$ is a minimal counterexample to $\d$-iterability. We then have that $\eta$ is a strong cutpoint of $\S$, $\S$ is $\eta$-sound and $\rho_\omega(\S)\leq \eta$.\\

\textbf{Assume that $\S$ has a Woodin cardinal $>\eta$.} Let $\nu$ be its least Woodin cardinal $>\eta$. Let $\P$ be the output of the fully backgrounded construction of $\M|\d[g]$ done over $\S|\eta$ using extenders with critical points $>\eta$. We now compare $\S|\nu$ with the construction producing $\P$. The $\P$-side of such a comparison doesn't move (for instance, see Lemma 2.11 of \cite{ATHM}). However, since $\S$ is not fully iterable, we need to describe a strategy for picking branches on the $\S$-side. Let $(\P^*_\xi, \P_\xi, E_\xi: \xi<\d)$ be the models of the aforementioned construction. 

Suppose then $\U\in \M|\d[g]$ is a tree of limit length that has been built on $\S$ via the aforementioned comparison process. We would like to describe a branch for it. As an inductive hypothesis, we maintain that $\VT^\frown \U$ is $\Lambda$-correct. Thus, the branch $b$ we pick for $\U$ has to have the property that $\VT^\frown \U^\frown\{b\}$ is $\Lambda$-correct. There can be at most one such branch. It is then enough to show that there is such a branch. The description of $b$ splits into two cases. 

First recall the definition of a fatal drop \cite[Definition 1.27]{ATHM}. Given a tree $\W$ on a premouse $\Q$ we say $\W$ has a fatal drop if there is $\a<lh(\W)$, $\xi$ and $\K\inseg \M^\T_\a$ such that $\xi$ is a strong cutpoint of $\K$, $\rho_\omega(\K)\leq \xi$ and $\T_{\geq \a}$ is an iteration of $\K$ above $\xi$. \\

\textbf{Case 1: $\U$ doesn't have a fatal drop}.\\

We have that there is some $\xi<\d$ such that $C(\U)\inseg \P_\xi$. Because $\d(\U)<\d$, we have that $\M\models ``\d(\U)$ is not a Woodin cardinal". It follows that there is a mouse $\Q$ over $C(\U)$ that is obtained via the $S$-construction that translates $\M$ into a mouse over $\P_\xi|\d(\U)$ such that $\Q$ is $\d(\U)$-sound, $\rho_\omega(\Q)\leq \d(\U)$ and $rud(\Q)\models ``\d(\U)$ is not a Woodin cardinal".  We claim that \\

\textit{Claim 3.} there is a branch $b$ of $\U$ such that $\Q(b, \U)$ exists and $\Q(b, \U)=\Q$. \\\\
\begin{proof}
Indeed, let $X\prec \M|(\d^+)^\M[g]$ be countable such that $\VT, \S, \U, \P_\xi , \Q \in X$, and letting $\pi_X:\N_X\rightarrow  \M|(\d^+)^\M[g]$ be the transitive collapse of $X$, $\pi^{-1}_X(\VT^\frown \U)$ is $\Lambda$-correct. Set $\bar{\Q}=\pi^{-1}(\Q)$. Notice that it follows from Claim 1 and \rprop{w is ub} that $\bar{\Q}$ is $\d$-iterable in $\N$. Let then $c=\Lambda(\pi_X^{-1}(\VT^\frown \U))$. We must have that $\Q(c, \bar{\U})$ exists and $\Q(c, \bar{\U})=\bar{\Q}$. By absoluteness $c\in \N_X$. It is now not hard to check that $b=_{def}\pi_X(c)$ is as desired. $b$ is the unique branch of $\U$ such that $\Q(b, \U)$ exists and $\Q(b, \U)=\Q$.\\

%
%
\end{proof}

\textbf{Case 2: $\U$ has a fatal drop.}\\

 Let $\xi<lh(\U)$ be such that the fatal drop happens at $\M^\U_\xi$. Let $\zeta$ and $\W\inseg \M^\U_\xi$ be such that $\M^\U_\xi|\zeta\inseg \W$, $\rho_\omega(\W)=\zeta$, $\zeta$ is a strong cutpoint of $\W$ and $\U_{\geq \xi}$ is an iteration tree on $\W$ above $\zeta$. Because $\S$ is a minimal counterexample to $\d$-iterability, we have that $\W$ is $\d$-iterable in $\M|\d[g]$. Let then $b$ be the branch of $\U$ according to the unique strategy of $\W$. Again a Skolem hull argument and Claim 1 show that $\VT^\frown \U^\frown \{ b\}$ is $\Lambda$-correct. \\
 
 This finishes our description of branches that payer $II$ plays in the comparison game between $\S|\nu$ and the construction producing $\P$. Let then $\U$ be the tree on $\S|\nu$ of maximal length constructed in the manner described above. 

Notice that for unboundedly many $\theta< \d$, $\P$ computes $\theta^+$ correctly. This is because $\X\in \P$ and if $\P^*$ is the output of the fully backgrounded construction of $\P$ done over $\X$ with large enough critical points then $\P^*$ computes unboundedly many successors correctly (this is a consequence of the internal covering property). 

It now follows that that  $lh(\U)\not=\d$ as then $C(\U)=\P$. Also, $\U$ must have a last model $\S^*$. Indeed, if $\U$ doesn't have a last model then it is of limit length. Because $\d(\U)<\d$, we have that $\M\models ``\d(\U)$ is not a Woodin cardinal", implying that our method of picking branches of $\U$ does produce a branch for $\U$. Because the $\S|\nu$-side lost the comparison, $\pi^\U$ must exist. Let  $\zeta=\pi^\U(\nu)$ (because $\VT^\frown \U$ is $\Lambda$-correct, $\U$ can be applied to $\S$). 

Because $\M\models ``\zeta$ is not a Woodin cardinal", we can find sound $\S^*|\zeta$-mouse $\W\in \M|\d[g]$ such that $\rho_\omega(\W)\leq \zeta$, $\W\models ``\zeta$ is a Woodin cardinal" and $rud(\W)\models ``\zeta$ is not a Woodin cardinal". Because we are in the no-fatal drop case, $\W=\S^*$.

%
%

We claim that in $\M[g]$, $\W$ is $\d$-iterable above $\zeta$.  Notice that because $\W$ is obtained via $S$-constructions, for some $\b$, $\W$ is the result of the translation of $\M||\b$ into a mouse over $\S^*|\zeta$ via $S$-constructions.  It follows that $\rho_{\omega}(\M|\b)\leq \zeta$. Since $\M||\b$ is $\d$-iterable above $\zeta$, $\W$ is also $\d$-iterable above $\zeta$ (see \rprop{everything is ub}). It follows from $S$-reconstructibility that for any $\xi\in (\zeta, \d)$, $\W$ can also be build by the fully backgrounded construction of $\M|\d[g]$ that uses extenders with critical points $>\xi$. Let now $\xi_0<\xi_1$ be $\M$-successor cardinals such that $\S^*|\zeta$ can be build by the fully backgrounded construction of $\M|\xi_0[g]$ that uses extenders with critical points $>\eta$, and $\W$ can be build by the fully backgrounded construction of $\M|\xi_1[g]$ that uses extenders with critical points $>\xi_0$. It follows that iterability of $\W$ can be reduced to the iterability of $\M|\xi_1$ for non-dropping trees that are above $\eta$. \rprop{everything is ub} then implies that $\W$ is $\d$-iterable in $\M[g]$.

This is a contradiction as  in $\M[g]$, $\S$ is not $\d$-iterable above $\eta$, while because $\pi^\U:\S\rightarrow \W$, we get that in $\M[g]$, $\S$ is in fact $\d$-iterable above $\eta$.

The case when $\S$ has no Woodin cardinals is very similar. Now we compare $\S$ with the fully backgrounded constructions producing a tree $\U$ on $\S$ such that $\VT^\frown \U$ is $\Lambda$-correct. Because $\S$ has no Woodin cardinals, handling limit stages of the construction of $\U$ is very similar. Assuming $\U$ has been built up to stage $\gg$, we consider, as above, two cases. If $\U\rest \gg$ has no fatal drops then we proceed as in the ``no fatal drop case" of the above argument. Otherwise, we proceed in the ``fatal drop" case of the above argument. We leave the details to the reader.   
\end{proof}

We believe that the project of characterizing in mice the exact cardinals $\k$ that permit stationary tower like embeddings with critical point $\k$ is a very nice project.

\bibliographystyle{plain}
\bibliography{Prepideals}

\end{document}